\documentclass[11pt]{amsart}
\usepackage[margin=10pt,font=small,labelfont=normalfont]{caption}
\usepackage{subcaption}

\usepackage{mathtools}

\usepackage[demo]{graphicx}

\usepackage{xcolor}

\usepackage[margin=1.1in]{geometry} 

\usepackage{graphics,graphicx}
\usepackage{pstricks,pst-node,pst-tree}

\usepackage{xcolor}
\theoremstyle{plain}

\usepackage{amsfonts}
\usepackage{amssymb}
\usepackage{amsthm}
\usepackage{amsmath}
\usepackage{enumerate}
\usepackage{hyperref}
\hypersetup{
    colorlinks=false,
    pdfborder={0 0 0},
}
\usepackage{amsrefs}
\usepackage{graphicx}
\usepackage{verbatim}
\usepackage{amscd}
\usepackage{mathtools}

\usepackage{arydshln}

\theoremstyle{plain}

\def\sqr#1#2{{\,\vcenter{\vbox{\hrule height.#2pt\hbox{\vrule width.#2pt
height#1pt \kern#1pt\vrule width.#2pt}\hrule height.#2pt}}\,}}

\newtheorem{proposition}{Proposition}

\newtheorem{theorem}[proposition]{Theorem}

\theoremstyle{definition}
\newtheorem{definition}[proposition]{Definition}
\theoremstyle{remark}

\newtheorem*{convention*}{Convention}

\newtheorem{remark}[proposition]{Remark}

\numberwithin{equation}{section}

\newcommand{\calC}{\hbox{$\mathcal C$}}
\newcommand{\calA}{\hbox{$\mathcal A$}}
\newcommand{\calH}{\hbox{$\mathcal H$}}

\newcommand{\calD}{\hbox{$\mathcal D$}}

\numberwithin{equation}{section}

\begin{document}

\title{Corrigendum to \\ Ternary operator categories}


\author[Pluta]{Robert Pluta}
\email{plutar@tcd.ie}
\address{Department of Mathematics, College of the Holy Cross, Worcester, MA, 01610, USA}

\author[Russo]{Bernard Russo}
\email{brusso@math.uci.edu}
\address{Department of Mathematics, UC Irvine, Irvine CA, USA}


\date{May 20, 2023}

\maketitle

\begin{abstract}
 We acknowledge an oversight in \cite{PluRusJMAA} and adjust the one section of that paper which is affected by it. 
 \end{abstract}


One of the results of \cite{PluRusJMAA}, namely,  \cite[Prop.\  2.7]{PluRusJMAA}, is not accurately formulated as stated. It is not true that the normed standard embedding ${\mathcal A}(M)$ of a C*-ternary ring $M$  is always a  C*-algebra. 
The correct version of Proposition 2.7 in \cite{PluRusJMAA}  would state that ${\mathcal A}(M)$ is the direct sum of a C*-algebra (corresponding to the positive sub-C*-ternary ring $M_+$ of $M$) and a Banach algebra ${\mathcal B}$ (corresponding to the negative sub-C*-ternary ring $M_-$ of $M$).   This will be made precise in a forthcoming paper \cite{PluRusAnti}. Moreover, the above mentioned oversight affects only the results in  \cite[Section 4]{PluRusJMAA}. The corresponding adjustments needed in \cite[Section 4]{PluRusJMAA} are detailed here in this corrigendum. In addition, we take this opportunity  to complete two proofs from \cite{PluRusJMAA} and provide an errata, so that all other results of \cite{PluRusJMAA} remain valid as stated.

The algebra ${\mathcal B}$ is a Banach algebra with a continuous involution and an approximate identity, and properties of this class of algebras is a topic worthy of further study.  We began that study in \cite{PluRusAnti} by showing, among other things, that it is semisimple.  
\medskip



\begin{center}{\bf Adjustments to  \cite[Section 4]{PluRusJMAA}}\end{center}

In what follows, we shall use some notation and some results from  \cite{PluRusJMAA}, making precise references to  \cite{PluRusJMAA} when necessary.
We first modify, as suggested in \cite{PluRusJMAA}, 
the definition of linking category. Recall that the morphism sets $(X,Y)_{\mathcal C}$ in a linear ternary category are associative triple systems.

\begin{definition}
(Modification of \cite[Def. 3.8]{PluRusJMAA})
Given a linear ternary category $\calC$, the {\it linking category} $A_{\calC}$  of $\calC$ is as follows. The objects of the category $A_{\calC}$ are the same as the objects of $\calC$.  The morphism set  $\operatorname{Hom}(X,Y)\ (=(X,Y)_{A_{\mathcal C}})$ is defined to be $\calA(X,Y)\  (={\mathcal A}((X,Y)_{\mathcal C})$, with composition  as follows. If $a\in \operatorname{Hom}(X,Y)$ and $b\in \operatorname{Hom}(Y,Z)$, then
$b\circ a$ must be $0$ unless $X=Y=Z$, in which case $b\circ a$ is defined to be the product $ab$ in $\calA(X,X)$. 
\end{definition}

\begin{remark}\label{rem:0420231}
The category $A_{\calC}$ can be considered as a ternary category (\cite[Definition 3.1]{PluRusJMAA}) under the composition $[abc]=ab^\#c$ and, by \cite[Lemma 1.3]{PluRusJMAA},  we obtain an injective linear ternary functor $F$ from $\calC$ to $A_{\calC}$ by associating the object $X$ of $\calC$ to the object $X$ of $A_{\calC}$ and the morphism $f\in (X,Y)$ to  the morphism $\left[
\begin{matrix}
0&f\\
0&0
\end{matrix}\right]\in\calA(X,Y).
$
\end{remark}

Recall that the morphism sets in a T*-category (\cite[Definition 4.6]{PluRusJMAA}) are C*-ternary rings (\cite[Section 2]{PluRusJMAA}), so they satisfy Zettl's representation theorem \cite[Theorem 2.1]{PluRusJMAA}.

\begin{definition}\label{def:0212231}
Given a T*-category $\calC$, the {\it positive (resp. negative) linking category} $A_{\calC}^{\pm}$  of $\calC$ is as follows. The objects of the category $A_{\calC}^{\pm}$ are the same as the objects of $\calC$.  The morphism set  $\operatorname{Hom}(X,Y)$ is defined to be $\operatorname{Hom}(X,Y)=\calA((X,Y)_{\pm})$, with composition  as follows. If $a\in \operatorname{Hom}(X,Y)$ and $b\in \operatorname{Hom}(Y,Z)$, then
$b\circ a$ must be $0$ unless $X=Y=Z$, in which case $b\circ a$ is defined to be the product $ab$ in $\calA((X,X)_{\pm})$. 
\end{definition}

Thus the  linking categories $A_{\calC}^{\pm}$ of the T*-category $\mathcal C$
  are subcategories of the linking category $A_{\calC}$  of the linear ternary category $\mathcal C$, and  $
  A_{\calC}=A_{\calC}^+\oplus A_{\calC}^-,
    $
    by \cite[Remark 4.13(ii)]{PluRusJMAA}.

\begin{theorem}\label{thm:0217231}
(Replacement for \cite[Theorem 4.11]{PluRusJMAA})
If $\calC$ is a T*-category  then $A_{\calC}^+$ is a C*-category (\cite[Definition 4.1]{PluRusJMAA}) and there is a faithful T*-functor $F$  from $\calC_+$ to $A_{\calC}^+$, the latter considered as a T*-category. 
\end{theorem}
\begin{proof}
It is clear that $A_{\calC}^+$, as defined in Definition~\ref{def:0212231}, is a linear non-unital category which, when considered as a ternary category, satisfies (ii), (iii) and (vi) in 
\cite[Def. 4.1]{PluRusJMAA}.
   Items (i), (iv), and (v)  of that definition
are tantamount  to the morphism sets of  $A_{\calC}^+$, namely, $\calA((X,Y)_+)$,  being normed as  C*-algebras. This fact is immediate from \cite[Theorem 2.7]{PluRusAnti}, and the faithful functor is the restriction of the functor defined in Remark~\ref{rem:0420231}.  
\end{proof}

Let $\rho=(\rho_0, \{\rho_{X,Y}\})$ be a linear ternary functor from a linear ternary category $\calC$ to a a linear ternary category $\calD$. We write $\rho=(\rho_0, \{\rho_{X,Y}:X,Y\hbox{ objects  of } \calC\})$, where $\rho_0$ maps objects of $\calC$ to objects of $\calD$, and $\rho_{X,Y}$ is a linear transformation from $(X,Y)_{\calC}$ to $(\rho_0(X),\rho_0(Y))_{\calD}$ satisfying 
\[
\rho_{X,W}(h\circ g^*\circ f)=\rho_{Z,W}\circ\rho_{Z,Y}(g)^*\circ \rho_{X,Y}(f),
\]
 for $f\in (X,Y)_{\calC},\  g\in (Z,Y)_{\calC},\hbox{ and } h\in (Z,W)_{\calC}$.
 
 If $\calC$ is a T*-category, 
 T*-subcategories $\calC_{\pm}$ are defined as follows. By Zettl's representation theorem, we have $
(X,Y)=(X,Y)_+\oplus (X,Y)_-
$  for each pair of objects $X,Y$ of $\calC$. The objects of 
$\calC_{\pm}$ are the same as the objects of $\calC$, and for such objects $X,Y$,
$$(X,Y)_{\calC_{\pm}}:=(X,Y)_{\pm}.$$
 It is clear that $\calC$ is isomorphic to $\calC_+\oplus \calC_-$ and that $A_{\calC}$ is isomorphic to
 $\calA_{C_+}\oplus \calA_{C_-}$ (cf.\  \cite[Remark 4.13(ii)]{PluRusJMAA}).

\begin{remark}
(Replacement for \cite[Remark 4.12]{PluRusJMAA})
If $\rho_0$ is injective,  then there is a linear functor $\widehat\rho$  from $A_{\calC}$ to $A_{\calD}$ which extends $\rho$. In particular, if $\calC$ and $\calD$ are T*-categories, then every T*-functor  from $\calC_+$ to $\calD_+$, which is injective on objects, extends to a C*-functor
  from $A_{\calC}^+$ to $A_{\calD}^+$.
\end{remark}
\begin{proof}
Recall that $(X,Y)_{\calC}$ is an associative triple system, and $\rho_{X,Y}$ is a homomorphism of 
$(X,Y)_{\calC}$ to $(\rho_0(X),\rho_0(Y))_{\calD}$, so by \cite[Lemma 2.6]{PluRusJMAA}, it extends to a *-homomorphism $\widehat\rho_{X,Y}={\mathcal A}(\rho_{X,Y})$ from $(X,Y)_{A_{\calC}}$  to 
$(\rho_0(X),\rho_0(Y))_{A_{\calD}} $.

Thus $\widehat\rho=(\widehat\rho_0, \{\widehat\rho_{X,Y}: X,Y\hbox{ objects of }A_{\calC}\})$, where $\widehat\rho_0(X)=\rho_0(X)$,  is the desired linear  functor from $A_{\calC}$ to $A_{\calD}$ whose restriction to ${\mathcal A}_{\mathcal C}^+$ is a C*-functor to ${\mathcal A}_{\mathcal D}^+$.
\end{proof}

\begin{theorem}
(Replacement for \cite[Theorem 4.14]{PluRusJMAA})
Let $\calC$ be a T*-category.  Then there is  a faithful T*-functor  from $\calC_+$ to the T*-category $\calH$ of Hilbert spaces and continuous linear maps (\cite[Example 4.2, Remark 4.8]{PluRusJMAA}).
\end{theorem}
\begin{proof} By \cite[Theorem 4.5]{PluRusJMAA}, there is a faithful C*-functor $G_{+}$ from $A_{\calC_{+}}$ to $\calH$.  With $A_{\calC_{+}}$ considered as  a T*-category, we have that $G_{+}$ is a T*-functor from $A_{\calC_{+}}$ to $\calH$. By Theorem~\ref{thm:0217231}, there is a faithful T*-functor $F_{+}$ from $\calC_{+}$  to $A_{\calC_{+}}$, and it suffices to consider $ H=G_+\circ F_+$.
\end{proof}

Let $V$ be a C*-ternary ring. By \cite[Theorem 2.7]{PluRusAnti}, $V$ is the off-diagonal corner of a Banach algebra  with continuous involution (a C*-algebra  if $V$ is a TRO),  $\calA(V)$, where
\[
\calA(V)=\left[\begin{matrix}
L&V\\
\overline{V}&R
\end{matrix}\right] ,
\]  
and $L=L(V)$ and $R=R(V)$ are C*-algebras.  Consider
\begin{equation}\label{eq:0915201}
\tilde A(V)=\left[\begin{matrix}
M(L)&V\\
\overline{V}&M(R)
\end{matrix}\right],
\end{equation} 
where $M(L)$ and $M(R)$ are the multiplier algebras of $L$ and of $R$.
Zettl has shown in \cite[Prop. 4.9]{Zettl83} that if $V$ is a W*-ternary ring  (\cite[Section 2]{PluRusJMAA}), then $M(R)$ and consequently $M(L)$ are W*-algebras.

Recall that a T*-category is a {\it TW*-category} if each morphism set is a dual Banach space, and that for objects $X,Y$ in a TW*-category, $(X,Y)$ is a W*-ternary ring.

\begin{definition}
\label{def:0219231}
(Replacement for \cite[Def. 4.22]{PluRusJMAA})
Given a TW*-category $\calC$ (\cite[Definition 4.6]{PluRusJMAA}), the {\it linking W*-category} $\tilde A_{\calC}^+$  of $\calC$ is as follows. The objects of the category $\tilde A_{\calC}^+$ are the same as the objects of $\calC$.  The morphism set  $\operatorname{Hom}(X,Y)$ is defined to be $\tilde \calA((X,Y)_+)$, as in (\ref{eq:0915201}) with $V=(X,Y)_+$, and with composition  as follows. If $a\in \operatorname{Hom}(X,Y)$ and $b\in \operatorname{Hom}(Y,Z)$, then
$b\circ a$ must be $0$ unless $X=Y=Z$, in which case $b\circ a$ is defined to be the product $ab$ in $\tilde \calA((X,X)_+)$. 
\end{definition}

\begin{theorem}\label{thm:0831202}
(Replacement for \cite[Theorem 4.23]{PluRusJMAA})
If $\calC$ is a TW*-category  then $\tilde A_{\calC}^+$ is a W*-category (C*-category with morphism sets being dual spaces) and there is a faithful TW*-functor $F$  from $\calC_+$ to $\tilde A_{\calC}^+$, the latter considered as a TW*-category. 
\end{theorem}
\begin{proof} 
It is clear that $\tilde A_{\calC}^+$, as defined in Definition~\ref{def:0219231}, is a linear non-unital category which, when considered as a ternary category, satisfies (ii), (iii) and (vi) in 
\cite[Def. 4.6]{PluRusJMAA}.
   Items (i), (iv), and (v)  of that definition are tantamount  to the morphism sets of  $\tilde A_{\calC}^+$, namely, $ \calA((X,Y)_+)$,  being normed as  W*-algebras. This fact is immediate from \cite[Prop. 4.21]{PluRusJMAA} (see Errata 7).
\end{proof}

\begin{theorem}
(Replacement for  \cite[Theorem 4.25]{PluRusJMAA})
Let $\calC$ be a TW*-category.  Then there is  a faithful TW*-functor  from $\calC_+$ to the TW*-category $\calH$ of Hilbert spaces and continuous linear maps.
\end{theorem}
\begin{proof} By \cite[Prop. 2.13]{GLR85}, there is a faithful W*-functor $F$ from $\tilde A_{\calC}^+$ to $\calH$.  With $\tilde A_{\calC}^+$ considered as  a TW*-category, we have that $G$ is a TW*-functor from $\tilde A_{\calC}^+$ to $\calH$. By Theorem~\ref{thm:0831202}, there is a faithful TW*-functor $F$ from $\calC^+$  to $\tilde A_{\calC}^+$, and it suffices to consider $ H=G\circ F$.
\end{proof}

\begin{center}{\bf Two completed proofs for \cite{PluRusJMAA}}\end{center}

\begin{remark}
(Completion of  the proof of  \cite[Rem. 4.10]{PluRusJMAA})
In the proof of  \cite[Rem. 4.10]{PluRusJMAA},  only items (i)-(iv) of the five conditions in \cite[Def. 4.6]{PluRusJMAA} were proved. The remaining item  (v) (as well as item (iv)) follow immediately from \cite[Prop.\ 4.2]{PluRusAnti}.
\end{remark}

In the proof of \cite[Prop.\ 5.13]{PluRusJMAA}, it was stated incorrectly that C*-ternary rings are JB*-triples, and therefore, since JB*-triples satisfy Pelczynski's property V, so do C*-ternary rings. However, as pointed out in \cite[Theorem 2.8]{PluRusAnti}, a C*-ternary ring $M=M_+\oplus M_-$ is (isomorphic to) a JB*-triple if and only if $M_-=0$. Thus a C*-ternary ring is a JB*-triple if and only if it is isomorphic to a TRO. Using this fact, we can prove the following proposition, which completes the proof of  \cite[Prop.\ 5.13]{PluRusJMAA}.

\begin{proposition}\label{prop:0420231}
A C*-ternary ring $(M,[\cdot,\cdot,\cdot])$ satisfies Pelczynski's property V.
\end{proposition}
\begin{proof}
According to the theorem of Zettl  \cite[Theorem 2.1]{PluRusJMAA}, there is a bounded operator $T$ on $M$ such that 
$(M,T\circ [\cdot,\cdot,\cdot])$ is isomorphic to a TRO. Therefore $(M,T\circ [\cdot,\cdot,\cdot])$ satisfies Pelczynski's property V. Since $(M,[\cdot,\cdot,\cdot])$ and $(M,T\circ [\cdot,\cdot,\cdot])$ are identical as Banach spaces, the proposition is proved.
\end{proof}

\begin{center}{\bf Errata for \cite{PluRusJMAA}}\label{sec:0410235}\end{center}
\begin{enumerate}
\item In \cite[Lemma 1.1(i)]{PluRusJMAA}, there is an extra left parenthesis.
\item In \cite[Lemma 1.2(i)]{PluRusJMAA}, ``$(L,R)$-bimodule'' should read ``$(L,R^{op})$-bimodule,'' and in \cite[Lemma 1.2(ii)]{PluRusJMAA}
``$(R,L)$-bimodule'' should read ``$(R^{op},L)$-bimodule.''
\item In the last sentence of the penultimate paragraph preceding \cite[Prop.\ 2.3]{PluRusJMAA}, ``$(L,R)$-bimodule'' should read 
``$(L,R^{op})$-bimodule''  and ``$(\frak B,\frak A)$-bimodule'' should read ``$(\frak B,{\frak A}^{op})$-bimodule.''
\item In the first line of the proof of \cite[Prop.\ 2.4]{PluRusJMAA}, ``$r(f,g)$'' should read  ``$r(g,f)$.'' 
\item There are several misprints in \cite[Example 4.16]{PluRusJMAA}, which are corrected in \cite[Prop.\ 2.2]{PluRusAnti}.
\item In the paragraph following \cite[Prop.\ 4.20]{PluRusJMAA}, $V$ is the off-diagonal corner of a Banach algebra, not necessarily a C*-algebra, unless $V$ is a TRO.
\item In \cite[Prop.\ 4.21]{PluRusJMAA}, $V$ is the off-diagonal corner of a dual Banach algebra, not necessarily a W*-algebra, unless $V$ is a W*-TRO.
\item In the last sentence of the proof of \cite[Prop.\ 4.21]{PluRusJMAA}, $\widetilde{A}(V)$ is a dual Banach algebra, not necessarily a W*-algebra, unless $V$ is a W*-TRO.
\end{enumerate}


\begin{bibdiv}
\begin{biblist}

\bib{GLR85}{article}{
   author={Ghez, P.},
   author={Lima, R.},
   author={Roberts, J. E.},
   title={$W^\ast$-categories},
   journal={Pacific J. Math.},
   volume={120},
   date={1985},
   number={1},
   pages={79--109},
   issn={0030-8730},
   review={\MR{808930}},
}





\bib{PluRusJMAA}{article}{
   author={Pluta, Robert},
   author={Russo, Bernard},
   title={Ternary operator categories},
   journal={J. Math. Anal. Appl.},
   volume={505},
   date={2022},
   number={2},
   pages={Paper No. 125590, 37},
   issn={0022-247X},
   review={\MR{4304265}},
   doi={10.1016/j.jmaa.2021.125590},
}


\bib{PluRusAnti}{article}{
   author={Pluta, Robert},
   author={Russo, Bernard},
   title={Anti-C*-algebras},
   journal={preprint},
}


\bib{Zettl83}{article}{
   author={Zettl, Heinrich},
   title={A characterization of ternary rings of operators},
   journal={Advances in Math.},
   volume={48},
   number={2},
   date={1983},
   pages={117--143},
  }


 
\end{biblist}
\end{bibdiv}
\end{document}